\documentclass[10pt,thmsa]{article}
\usepackage{amsmath, latexsym, amsfonts, amssymb, amsthm, amscd, verbatim}

\newtheorem{theorem}{Theorem}

\newtheorem{corollary}[theorem]{Corollary}

\newtheorem{definition}[theorem]{Definition}

\newtheorem{lemma}[theorem]{Lemma}

\newtheorem{proposition}[theorem]{Proposition}
\newtheorem{remark}[theorem]{Remark}

\newcommand{\D}{{\rm d}}
\newcommand{\E}{{\rm E}}

\newcommand{\p}{\Bbb{P}}

\newcommand{\e}{\Bbb{E}}

\def\QED{\hfill\vrule height 1.5ex width 1.4ex depth -.1ex \vskip20pt}

\begin{document}

\title{Refracted L\'evy processes.}

\maketitle

\begin{center}
 {\large  A.E. Kyprianou$^{*,}$\footnote{\noindent$^{,2}$ Department of Mathematical Science, University of Bath. {\sc Bath, BA2 7AY. United Kingdom.}
$^{1}$E-mail: a.kyprianou@bath.ac.uk, \,$^2$E-mail: rll22@bath.ac.uk.\\
$^*$ Corresponding author.} and R.L. Loeffen$^2$}
\end{center}
\vspace{0.2in}

\begin{abstract}
Motivated by classical considerations from risk theory, we investigate boundary crossing problems for refracted L\'evy processes. The latter is a L\'evy process whose dynamics change by subtracting off a fixed linear drift (of suitable size) whenever the aggregate process is above a pre-specified level. More formally, whenever it exists, a refracted L\'evy process is described by the unique strong solution to the stochastic differential equation
\[
\D U_t = - \delta \mathbf{1}_{\{U_t >b\}}\D t + \D X_t
\]  
where $X=\{X_t :t\geq 0\}$ is a L\'evy process with law $\mathbb{P}$ 
and $b, \delta\in \mathbb{R}$ such that the resulting process $U$ may visit the half line $(b,\infty)$ with positive probability. 
We consider in particular the case that $X$ is spectrally negative and establish a suite of identities for the case of one and two sided exit problems. All identities can be written in terms of the $q$-scale function of the driving L\'evy process and its perturbed version describing motion above the level $b$. We remark on a number of applications of the obtained identities to (controlled) insurance risk processes.


\bigskip


\end{abstract}

\vspace{0.5cm}


\section{Introduction.}

In this paper we are interested in understanding the dynamics of a one-dimensional L\'evy 
 process when its path is perturbed in a simple way. Informally speaking, a linear drift at rate $\delta>0$ is subtracted from the increments of a L\'evy process whenever it exceeds a pre-specified positive level. 
 More formally, suppose that $X=\{X_t : t\geq 0\}$ is L\'evy process. If we denote the level by $b>0$, a natural way to model such processes is to consider them as solutions to the stochastic differential equation
\begin{equation}
U_t = X_t - \delta \int_0^t \mathbf{1}_{\{U_s >b\}}\D s, \, t\geq 0
\label{SDE}
\end{equation}
assuming that at least a unique weak solution exists and  such that $U=\{U_t : t\geq 0\}$ visits $(b,\infty)$ with positive probability.

As a first  treatment of (\ref{SDE}) we shall restrict ourselves to the case that $X$ is a process with no positive jumps and such that $-X$ is not a subordinator (also henceforth referred to as spectrally negative L\'evy processes).  As a special case of the latter, suppose that $X$ may be written in the form 
\begin{equation}
X_t = c t- S_t, \, t\geq 0
\label{sn}
\end{equation}
 where $c>0$ is a constant and $S=\{S_t: t\geq 0\}$ is a pure jump compound Poisson subordinator.
In that case it is easy to see that, under the hypothesis
$c>\delta$
a solution to (\ref{SDE})
may be constructed pathwise utilizing the the fact that $b$ is always crossed by $X$ from below on the path of a  linear part of the trajectory at a discrete set of times and is always crossed by $X_t-\delta t$ from above by a jump. 
Note that the trajectory of the process $U$ is piecewise linear and `bent' as it crosses the level $b$ in the fashion that a light ray refracts from one medium to another. Inspired by this mental picture, 
we refer to solutions of (\ref{SDE}) when the driving process $X$ is a general one dimensional L\'evy process as a {\it refracted L\'evy process}.\footnote{See for example the diagram on p80 of \cite{GS2006} and the text above it which also makes reference to `refraction' in the case of compound Poisson jumps. The article \cite{GS2006BM} also uses the terminology `refraction' for the case that $X$ is a linear Brownian motion.}

 The special case (\ref{sn}) with compound Poisson jumps described above may also be seen as an example of a Cram\'er-Lundberg process as soon as $\mathbb{E}(X_1)>0$. This provides a specific motivation for the study of the dynamics of (\ref{SDE}). Indeed very recent studies of problems related to ruin in insurance risk has seen some preference to working with general spectrally negative L\'evy processes in place of the classical Cram\'er-Lundberg process (which is itself an example of the former class). See for example \cite{APP2007,F1998, HPSV2004a,  HPSV2004b, KKM2004, KK2006, KP2007, RZ2007, SV2007}. This preference is largely thanks to  the robust mathematical theory which has been developed around certain path decompositions of such processes as well as the meaningful interpretation of the general spectrally negative L\'evy process as an insurance risk process (see for example the  discussion in Section \ref{applications} or \cite{KKM2004, SV2007}).  

Under such a general model, the solution to the stochastic differential equation (\ref{SDE}) may now be thought of the aggregate of the insurance risk process when dividends are paid out at a rate $\delta$ whenever it exceeds the level $b$. 
Quantities which have been of persistent interest in the literature invariably pertain to the behaviour of (\ref{SDE}) up to the ruin time $\kappa^-_0=\inf\{t>0 : U_t <0\}$. For example, 
the probability of ruin, $\mathbb{P}_x(\kappa^-_0<\infty)$, the net present value of the dividends paid out until ruin, 
$
\mathbb{E}_x\left(\int_0^{\kappa^-_0} e^{-qt}\delta\mathbf{1}_{\{U_t >b\} } \D s \right),
$
where $q>0$,
and the overshoot and undershoot at ruin, $\mathbb{P}_x(U_{\kappa^-_0}\in A,U_{\kappa^-_0-}\in B)$ where $A\subset(-\infty, 0)$, $B\subset [0,\infty)$ and $U_{\kappa^-_0-}=\lim_{t\uparrow\kappa^-_0}U_{t}$.
Whilst expressions for the expected discounted value of the dividends, the Laplace transform of the ruin probability and the joint law of the undershoot and overshoot have been established before for refracted L\'evy processes (cf. \cite{linpavlova}, \cite{wan2007}, \cite{zhangzhouguo}, \cite{zhou2004}, \cite{zhoudiscussion}) none of them go beyond the case of a compound Poisson jump structure. Moreover, existing identities in these cases are not often written in the modern language of scale functions (defined in Section \ref{main} below). The latter  has some advantage given the analytical properties and families of examples that are now known for such functions (cf. \cite{HK2007, KR2007}).

\bigskip

Our objectives in this paper are three fold. Firstly to show that refracted L\'evy processes exist as a unique solution to (\ref{SDE}) in  the  strong sense  whenever $X$ is a spectrally negative L\'evy process (establishing the existence and uniqueness turns out to be not as simple as (\ref{SDE}) looks for some cases of driving process $X$). Secondly to study  their dynamics by establishing a suite of identities, written in terms of scale functions, related to one and two sided exit problems  and thirdly to cite the relevance of such identities in context of a number of recent and classical applications of spectrally negative L\'evy processes within the context of ruin problems.

The remainder of the paper is structured as follows. In the next section we compile all of our main results together. Principally these consist of showing the existence and uniqueness of solutions to (\ref{SDE}) which turns out to be in the strong sense. The principle difficulties that arise in handling (\ref{SDE}) lie with the case that $X$ has unbounded variation paths with no Gaussian part which seemingly falls outside of many standard results on existence and uniqueness of solutions to stochastic differential equations driven by L\'evy processes. 
Then in  Sections 3-9 we give the proofs of our main results. Finally, in section \ref{applications} we return to the discussion on applications in (controlled) risk processes where explicit examples are given.

\section{Main results}\label{main}

Henceforth the process $(X,\mathbb{P})$  will always denote a spectrally negative L\'evy process.  It is well known that spectral negativity allows us to talk about the Laplace exponent $\psi(\theta)=\log\mathbb{E}(e^{\theta X_1})$ for $\theta\geq 0$. Further the Laplace exponent is known to necessarily take the form
\begin{equation}
\psi(\theta) =  \left\{\frac{1}{2}\sigma^2 \theta^2 \right\} + \left\{
\gamma\theta  - \int_{(1,\infty)}(1- e^{-\theta x})\Pi(\D x)
\right\} 
- \left\{\int_{(0,1)} (1-  e^{-\theta x} - \theta x) \Pi(\D x)\right\}
\label{generalCR}
\end{equation}
for $\gamma\in\mathbb{R}$,  $\sigma\geq 0$ and  L\'evy measure $\Pi$ satisfying
 \[
 \Pi(-\infty,0)=0 \text{ and } \int_{(0,\infty)}(1\wedge x^2)\Pi(\D x)<\infty
 \]
  (even though $X$ only has negative jumps, for convenience we choose the L\'evy measure to have only mass on the positive instead of the negative half line). 
Note that when $\Pi(0,\infty)=\infty$ the process $X$ enjoys a countably infinite number of jumps over each finite time horizon.
We shall also denote by $\{\mathbb{P}_x: x\in\mathbb{R}\}$  probabilities of $X$ such that under $\mathbb{P}_x$, the process $X$ is issued from $x$. Moreover, $\mathbb{E}_x$ will be the expectation operator associated to $\mathbb{P}_x$. For convenience in the case that $x=0$ we shall always write $\mathbb{P}$ and $\mathbb{E}$ instead of $\mathbb{P}_0$ and $\mathbb{E}_0$.

We need the following hypothesis which will be in force throughout the remainder of the paper:

\bigskip

{\bf (H)}  {\it the constant $0<\delta< \gamma + \int_{(0,1)}x\Pi(\D x)$ if $X$ has paths of bounded variation.}

\bigskip

\noindent Note that when $X$ is a spectrally negative L\'evy  process with bounded variation paths, it can always be written in the form (\ref{sn}) where $c>0$ and $S$ is a pure jump subordinator.
In that case, one sees that the hypothesis  (H) simply says that ${c}>\delta>0$. 
Write $\mathcal{S}$ for the space of spectrally negative L\'evy processes satisfying (H). As well as writing $X\in\mathcal{S}$, we shall also abuse our notation and write $(\gamma,\sigma, \Pi)\in\mathcal{S}$ if $(\gamma,\sigma,\Pi)$ is the triplet associated to $X$. Below, our first result concerns existence and uniqueness of solutions to (\ref{SDE}).

\begin{theorem}\label{rongsitu}
For a fixed $X_0=x\in\mathbb{R}$, there exists a unique strong solution  to (\ref{SDE}) within the class $\mathcal{S}$.
\end{theorem}

\begin{remark}\rm\label{strongsigma}
The existence of a unique strong solution to (\ref{SDE}) is, to some extent, no surprise within the class  of solutions driven by a general L\'evy proceses (not necessarily spectrally negative) with non-zero Gaussian component. Indeed for the latter class, existence of a strong unique solution is known,  for example, from the work of  Veretennikov \cite{ver} and Theorem 305 of the Monograph of Situ \cite{Situ2005}. 
 The strength of Theorem \ref{rongsitu} thus lies in dealing with the case that $X\in\mathcal{S}$ with no Gaussian component. In fact it will turn out that the real difficulties lie with the case that $X$ has paths of unbounded variation with no Gaussian part. Such stochastic differential equations, in particular with drift coefficients which are neither Lipschitzian nor continuous but just bounded and measurable, are called degenerate and less seems to be known about them in the literature for the case of a driving L\'evy process. See for example the remark proceeding Theorem III.2.34 on p159 of \cite{JS2003} as well as the presentation in \cite{Situ2005}. 
\end{remark}

\begin{remark}\rm
Standard arguments show that the existence of a unique strong solution to (\ref{SDE}) for each point of issue $x\in\mathbb{R}$, lead to the conclusion that $U$ is a Strong Markov Process. Indeed suppose that $T$ is a stopping time with respect to the natural filtration generated by $X$. Then  define a process $\widehat{U}$ whose dynamics are those of $\{U_t: t\leq T\}$ issued from $x$ and, on the event that $\{T<\infty\}$,   it evolves on the time horizon $[T,\infty)$ as the unique solution, say $\widetilde{U}$, to (\ref{SDE}) driven by the L\'evy process $\widetilde{X}= \{X_{T+s} - X_T : s\geq 0\}$ when issued from the random starting point $U_T$. Note that by construction, on $\{T<\infty\}$, the dependency of $\{\widehat{U}_t : t\geq T\}$ on $\{\widehat{U}_t : t\leq T\}$ occurs only through the value $\widehat{U}_T = U_T$. Note also that for $s>0$
\begin{eqnarray*}
 \widehat{U}_{T+t} &=&\widetilde{U}_t\\
 &=& \widehat{U}_{T} + \widetilde{X}_{t} -\delta \int_0^{t}\mathbf{1}_{\{\widetilde{U}_s >b\}} \D s\\
 &=&  x+ X_T -\delta \int_0^T \mathbf{1}_{\{U_s>b\}} \D s+ (X_{T+t} - X_T) -\delta \int_0^{t}\mathbf{1}_{\{\widehat{U}_{T+s} >b\}} \D s\\
 &=& x+ X_{T+t} -\delta\int_0^{T+t}\mathbf{1}_{\{\widehat{U}_s >b\}} \D s
\end{eqnarray*}
showing that $\widehat{U}$ solves (\ref{SDE}) issued from $x$. Given there is strong uniqueness of solutions to (\ref{SDE}), we may identify this solution to be $\widehat{U}$ and thus in possession of the Strong Markov Property.
\end{remark}

Before proceeding to the promised fluctuation identities we must first recall a few facts concerning {\it scale functions} for spectrally negative L\'evy processes, in terms of which all identities will be written.  For each $q\geq 0$ define $W^{(q)}:\mathbb{R}\rightarrow [0,\infty)$ such that $W^{(q)}(x)=0$ for all $x<0$ and  on $(0,\infty)$ $W^{(q)}$ is the unique continuous function with Laplace transform
\begin{equation}
\int_0^\infty e^{-\beta x}W^{(q)}(x)\,\D x = \frac{1}{\psi(\beta) - q}\label{Laplace}
\end{equation}
for all $\beta> \Phi(q)$, where $\Phi(q)=\sup\{\theta\geq 0 : \psi(\theta) = q\}$. For convenience, we write  $W$ instead of $W^{(0)}$.
Associated to the functions $W^{(q)}$ are the functions $Z^{(q)}:\mathbb{R}\rightarrow [1,\infty)$ defined by
\[
Z^{(q)}(x) = 1 + q \int_0^x W^{(q)}(y)\, \D y
\]
for $q\geq 0$.
Together, the functions $W^{(q)}$ and $Z^{(q)}$ are collectively known as scale functions and predominantly appear in almost all fluctuation identities
for spectrally negative L\'evy processes. Indeed, several such identities which are well known (cf. Chapter 8 of \cite{K2006}) are given in Theorem \ref{appendixthrm} in the Appendix and will be of repeated use throughout the remainder of the text.


Note also that by considering the Laplace transform of $W^{(q)}$, it is straightforward to deduce that $W^{(q)}(0+)=1/{c}$ when $X$ has bounded variation and therefore is (necessarily)  written in the form ${c} t  - S_t$ where $S=\{S_t : t\geq 0\}$ is a driftless subordinator and ${c}>0$. Otherwise $W^{(q)}(0+)=0$ for the case of unbounded variation. In all cases, if $X$ drifts to $\infty$ then $W(\infty) = 1/\mathbb{E}(X_1)$. 
In general the derivative of the scale function is well defined except for at most countably many points. However, when $X$ has unbounded variation  or $\Pi$ has no atoms, then for any $q\geq 0$, the restriction of $W^{(q)}$ to the positive half line belongs to $C^{1}(0,\infty)$.  See for example \cite{Lam} and \cite{KRS}. In \cite{CK} it was also shown that when $X$ has a Gaussian component ($\sigma>0$), then $W^{(q)}\in C^2(0,\infty)$.
Finally it is worth mentioning that as the Laplace exponent $\psi$ is continuous in its L\'evy triplet (continuity for the L\'evy measure is understood in the sense of weak convergence), it follows by the Continuity Theorem for Laplace transforms that $W^{(q)}$ is also continuous in its underlying L\'evy triplet. Moreover, performing an integration by parts, one obtains 
\[
 \int_{[0,\infty)} e^{-\beta x }W^{(q)}(\D x) = \frac{\beta}{\psi(\beta) - q}
\]
for all $\beta> \Phi(q)$ which, by the same reasoning as before, shows that $W^{(q)\prime}$ is also continuous in its underlying L\'evy triplet.

\bigskip

We are now ready to state our main conclusions with regard to certain fluctuation identities. 
 In all theorems, the process $U=\{U_t: t\geq 0\}$ is the solution to (\ref{SDE}) when driven by $X\in\mathcal{S}$ and the level $b>0$. We shall frequently refer to the stopping times
  \[
\kappa^+_a: = \inf\{t> 0 : U_t >a\} \text{ and }\kappa^-_0: = \inf\{t>0 : U_t < 0 \}.
\]
where $a>0$.
Further, let $Y=\{Y_t := X_t-\delta t: t\geq 0\}$. For each $q\geq 0$, $W^{(q)}$ and $Z^{(q)}$ are the $q$-scale functions associated with $X$ and $\mathbb{W}^{(q)}$ and $\mathbb{Z}^{(q)}$ is the $q$-scale function associated with $Y$. Moreover $\varphi$ is defined as the right inverse of the Laplace exponent of $Y$ so that
\[
 \varphi(q)  = \sup\{\theta \geq 0: \psi(\theta) - \delta\theta = q\}.
\]

\begin{theorem}[Two sided exit problem]\label{two-sided}
$\left. \right.$ 
\begin{description}
\item[(i)] For $q\geq0$ and $0\leq x,b\leq a$ we have 
\begin{equation}
\e_x(e^{-q\kappa^+_a}\mathbf{1}_{\{\kappa^+_a<\kappa^-_0\}})=\frac{W^{(q)}(x)+\delta \mathbf{1}_{\{x\geq b\}}\int_b^x\mathbb{W}^{(q)}(x-y)W^{(q)\prime}(y)\D y}{W^{(q)}(a)+\delta\int_b^a\mathbb{W}^{(q)}(a-y)W^{(q)\prime}(y)\D y}.
\label{useSkorokhod}
\end{equation} 
\item[(ii)]For $q\geq0$ and $0\leq x,b\leq a$ we have
\begin{multline*}
\mathbb{E}_x\left(\mathrm{e}^{-q\kappa^-_0}\mathbf{1}_{\{\kappa^-_0<\kappa^+_a\}}\right)=
Z^{(q)}(x)+\delta\mathbf{1}_{\{x\geq b\}}q\int_b^x \mathbb{W}^{(q)}(x-y)W^{(q)}(y)\mathrm{d}y \\
-\frac{Z^{(q)}(a)+\delta q\int_b^a\mathbb{W}^{(q)}(a-y)W^{(q)}(y)\mathrm{d}y}
{W^{(q)}(a)+\delta\int_b^a\mathbb{W}^{(q)}(a-y)W^{(q)\prime}(y)\mathrm{d}y}\\
\cdot
\left(W^{(q)}(x)+\delta\mathbf{1}_{\{x\geq b\}}\int_b^x \mathbb{W}^{(q)}(x-y)W^{(q)\prime}(y)\mathrm{d}y\right).
\end{multline*}
\end{description}
\end{theorem}

\begin{theorem}[One sided exit problem]\label{ruin} $\left. \right.$ 
\begin{description}
\item[(i)] For $q\geq0$ and $x,b\leq a$ we have
\[
\mathbb{E}_x(e^{-q\kappa^+_a }\mathbf{1}_{\{\kappa^+_a <\infty\}})=
\frac{\mathrm{e}^{\Phi(q) x}+\delta  \Phi(q) \mathbf{1}_{\{x\geq b\}} \int_{b}^{x}\mathrm{e}^{\Phi(q) z}\mathbb{W}^{(q)}(x-z)  \D z}
{\mathrm{e}^{\Phi(q) a}+  \delta\Phi(q)   \int_{b}^{a}\mathrm{e}^{\Phi(q) z} \mathbb{W}^{(q)}(a-z)  \D z}
\]

\item[(ii)] For $x,b\geq 0$ and $q>0$
\begin{multline*}
\mathbb{E}_x(e^{-q\kappa^-_0 }\mathbf{1}_{\{\kappa^-_0 <\infty\}}) = 
Z^{(q)}(x)+\delta \mathbf{1}_{\{x\geq b\}}q\int_b^x \mathbb{W}^{(q)}(x-y)W^{(q)}(y)\mathrm{d}y\\
 -\frac{q\int_b^\infty\mathrm{e}^{-\varphi(q)y}W^{(q)}(y)\mathrm{d}y}
{\int_b^\infty\mathrm{e}^{-\varphi(q)y}W^{(q)\prime}(y)\mathrm{d}y}
\left(W^{(q)}(x)+\delta\mathbf{1}_{\{x\geq b\}}\int_b^x \mathbb{W}^{(q)}(x-y)W^{(q)\prime}(y)\mathrm{d}y\right).
\end{multline*}
If in addition $0<\delta<\mathbb{E}(X_1)$, then letting $q\downarrow 0$ one has the ruin probability
\[
\p_x(\kappa^-_0<\infty) =1 - \frac{\mathbb{E}(X_1) - \delta }{1 - \delta W(b)}\left(W(x) + \delta\mathbf{1}_{\{x\geq b\}}\int_b^x \mathbb{W}(x-y) W'(y)\mathrm{d}y\right).
\]
\end{description}
\end{theorem}

\begin{theorem}[Resolvents]
\label{dividends}
Fix the Borel set $B\subseteq \mathbb{R}$.
\begin{description}
\item[(i)]
For $q\geq0$ and $0\leq x,b\leq a$, 
\begin{eqnarray}
\lefteqn{\mathbb{E}_x\left(\int_0^\infty e^{-qt}\mathbf{1}_{\{U_t \in B, \, t<\kappa^-_0\wedge \kappa^+_a\} } \D s \right)}&&\notag\\
&&= \int_{B\cap[b,a]} \left\{  \frac{W^{(q)}(x) + \delta \mathbf{1}_{\{x\geq b\}}\int_{b}^x\mathbb{W}^{(q)}(x-z)W^{(q)\prime}(z)\mathrm{d}z}{ W^{(q)}(a) + \delta\int_{b}^a\mathbb{W}^{(q)}(a-z)W^{(q)\prime}(z)\mathrm{d}z }   \mathbb{W}^{(q)}(a-y) -\mathbb{W}^{(q)}(x-y) \right\} \mathrm{d}y \notag\\
&&\hspace{0.5cm} + \int_{B\cap[0,b)} \Bigg\{   \frac{W^{(q)}(x) + \delta \mathbf{1}_{\{x\geq b\}}\int_{b}^x\mathbb{W}^{(q)}(x-z)W^{(q)\prime}(z)\mathrm{d}z}{ W^{(q)}(a) + \delta\int_{b}^a\mathbb{W}^{(q)}(a-z)W^{(q)\prime}(z)\mathrm{d}z }  \notag\\
&&\hspace{1.5cm} \cdot  \left( W^{(q)}(a-y) + \delta\int_{b}^{a}\mathbb{W}^{(q)}(a-z)W^{(q)\prime}(z-y)\mathrm{d}z \right) \notag\\
&&\hspace{2cm} - \left( W^{(q)}(x-y) + \delta\mathbf{1}_{\{x\geq b\}}\int_{b}^{x}\mathbb{W}^{(q)}(x-z)W^{(q)\prime}(z-y)\mathrm{d}z \right) \Bigg\} \mathrm{d}y .
\label{2sidedres}
\end{eqnarray}
\item[(ii)]For $x,b\geq 0$ and $q>0$
\begin{eqnarray*}
\lefteqn{\mathbb{E}_x\left(\int_0^\infty e^{-qt}\mathbf{1}_{\{U_t \in B, \, t<\kappa^-_0\} } \D s \right)}&&\\
&&=  \int_{B\cap[b,\infty)} \left\{   \frac{ W^{(q)}(x) + \delta\mathbf{1}_{\{x\geq b\}}\int_{b}^x\mathbb{W}^{(q)}(x-z)W^{(q)\prime}(z)\mathrm{d}z}{\delta \int_{b}^\infty\mathrm{e}^{-\varphi(q)z}W^{(q)\prime}(z)\mathrm{d}z}  \mathrm{e}^{-\varphi(q) y}  -\mathbb{W}^{(q)}(x-y) \right\} \mathrm{d}y \\
&& \hspace{0.5cm}+ \int_{B\cap[0,b)} \Bigg\{  \frac{\int_{b}^\infty\mathrm{e}^{-\varphi(q)z}W^{(q)\prime}(z-y)\mathrm{d}z }{ \int_{b}^\infty\mathrm{e}^{-\varphi(q)z}W^{(q)\prime}(z)\mathrm{d}z} \left( W^{(q)}(x) + \delta\mathbf{1}_{\{x\geq b\}}\int_{b}^x\mathbb{W}^{(q)}(x-z)W^{(q)\prime}(z)\mathrm{d}z \right)  \\
&& \hspace{1.5cm}- \left( W^{(q)}(x-y) + \delta\mathbf{1}_{\{x\geq b\}}\int_{b}^{x}\mathbb{W}^{(q)}(x-z)W^{(q)\prime}(z-y)\mathrm{d}z \right) \Bigg\} \mathrm{d}y.
\end{eqnarray*}

\item[(iii)] For $x,b\leq a$ and $q\geq0$
\begin{eqnarray*}
\lefteqn{\mathbb{E}_x\left(\int_0^\infty e^{-qt}\mathbf{1}_{\{U_t \in B, \, t<\kappa^+_a\} } \D s \right)}\\
&= &  \int_{B\cap[b,a]} \left\{  \frac{\mathrm{e}^{\Phi(q) x}+\delta  \Phi(q)  \mathbf{1}_{\{x\geq b\}} \int_{b}^{\infty}\mathrm{e}^{\Phi(q) z}\mathbb{W}^{(q)}(x-z)  \D z}
{\mathrm{e}^{\Phi(q) a}+  \delta\Phi(q)   \int_{b}^{\infty}\mathrm{e}^{\Phi(q) z} \mathbb{W}^{(q)}(a-z)  \D z}   \mathbb{W}^{(q)}(a-y) -\mathbb{W}^{(q)}(x-y) \right\} \mathrm{d}y \\
&& + \int_{B\cap(-\infty,b)} \Bigg\{   \frac{\mathrm{e}^{\Phi(q) x}+\delta  \Phi(q) \mathbf{1}_{\{x\geq b\}} \int_{b}^{\infty}\mathrm{e}^{\Phi(q) z}\mathbb{W}^{(q)}(x-z)  \D z}
{\mathrm{e}^{\Phi(q) a}+  \delta\Phi(q)   \int_{b}^{\infty}\mathrm{e}^{\Phi(q) z} \mathbb{W}^{(q)}(a-z)  \D z}  \\
&& \hspace{2cm}\cdot  \left( W^{(q)}(a-y) + \delta\int_{b}^{a}\mathbb{W}^{(q)}(a-z)W^{(q)\prime}(z-y)\mathrm{d}z \right) \\
&&\hspace{3cm} - \left( W^{(q)}(x-y) + \delta\mathbf{1}_{\{x\geq b\}}\int_{b}^{x}\mathbb{W}^{(q)}(x-z)W^{(q)\prime}(z-y)\mathrm{d}z \right) \Bigg\} \mathrm{d}y .
\end{eqnarray*}

\item[(iv)] For $x,b\in\mathbb{R}$ and $q> 0$, 
\begin{eqnarray*}
\lefteqn{\mathbb{E}_x\left(\int_0^\infty e^{-qt}\mathbf{1}_{\{U_t \in B\} } \D s \right)}
  && \\
  &&= \int_{B\cap[b,\infty)} \Bigg\{  \left( \mathrm{e}^{\Phi(q)(x-b)} + \delta\Phi(q)\mathrm{e}^{-\Phi(q)b}\mathbf{1}_{\{x\geq b\}}\int_{b}^{x}\mathrm{e}^{\Phi(q)z}\mathbb{W}^{(q)}(x-z)\mathrm{d}z \right)  \\
&&\hspace{3cm} \cdot    \frac{   \varphi(q)-\Phi(q) }{ \delta \Phi(q)  } 
\mathrm{e}^{-\varphi(q)(y-b)} - \mathbb{W}^{(q)}(x-y) \Bigg\} \mathrm{d}y \\
&& + \int_{B\cap(-\infty,b)} \Bigg\{  \left( \mathrm{e}^{\Phi(q)(x-b)} + \delta\Phi(q)\mathrm{e}^{-\Phi(q)b} \mathbf{1}_{\{x\geq b\}}\int_{b}^{x}\mathrm{e}^{\Phi(q)z}\mathbb{W}^{(q)}(x-z)\mathrm{d}z \right)  \\
&&   \hspace{3cm}\cdot \frac{   \varphi(q)-\Phi(q) }{ \Phi(q)  } 
\mathrm{e}^{\varphi(q)b}\int_{b}^\infty\mathrm{e}^{-\varphi(q)z}W^{(q)\prime}(z-y)\mathrm{d}z \\
&& \hspace{4cm}- \left( W^{(q)}(x-y) + \delta \mathbf{1}_{\{x\geq b\}}\int_{b}^{x}\mathbb{W}^{(q)}(x-z)W^{(q)\prime}(z-y)\mathrm{d}z \right)  \Bigg\} \mathrm{d}y.
\end{eqnarray*}

\end{description}
\end{theorem}

\begin{theorem}[Creeping]\label{creep}
For all $x,b\geq 0$ and $q>0$,
\begin{eqnarray*}
\lefteqn{\mathbb{E}_x\left( \mathrm{e}^{-q\kappa_0^-}\mathbf{1}_{\{ U_{\kappa_0^-}=0 \}}   \right) }\\
&&= \frac{\sigma^2}{2}  \Bigg\{    W^{(q)\prime}(x) + \delta \mathbf{1}_{\{x\geq b\}} \int_{b}^x\mathbb{W}^{(q)}(x-z)W^{(q)\prime\prime}(z)\mathrm{d}z  \\
&& \hspace{2cm} - \frac{ \int_b^\infty\mathrm{e}^{-\varphi(q)z}W^{(q)\prime\prime}(z)\mathrm{d}z }{ \int_b^\infty\mathrm{e}^{-\varphi(q)z}W^{(q)\prime}(z )\mathrm{d}z}    \left( W^{(q)}(x) + \delta \mathbf{1}_{\{x\geq b\}} \int_{b}^{x}\mathbb{W}^{(q)}(x-z)W^{(q)\prime}(z)\mathrm{d}z \right) \Bigg\}
\end{eqnarray*}
where the right hand side should be understood to be equal to zero when $\sigma=0$.  
\end{theorem}

\begin{remark}[Identities in Theorems \ref{dividends} and \ref{creep} when $q=0$]\label{limitingsense}\rm
In the previous two theorems the parameter $q$ was  taken to be strictly positive for some of the identities.
The case that $q=0$ can be handled by taking limits as $q=0$ on both left and right hand sides of these identities. 

\end{remark}

\begin{remark}\rm
Note that in the above (and subsequent) expressions the derivative of the scale function appears, despite the fact that in general $W^{(q)\prime}$ may not be well defined for a countable number of points. However, since $W^{(q)\prime}$ only appears in the integrand of an ordinary Lebesgue integral, this does not present a problem.
\end{remark}

\begin{remark}\rm
As is the case with any presentation which expresses identities in terms of scale functions of spectrally negative L\'evy processes, one may argue that one has only transferred the issue of `solving the problem' into finding explicit examples of scale functions. Although in general  scale functions are only semi-explicitly known through their Laplace transform, there are now quite a number of cases for which they can be calculated explicitly. See for example \cite{HK2007} and \cite{KR2007} for an updated account including  a variety of new, explicit  examples.  

For the cases where no explicit formula is known for the scale function, \cite{Rogers2000} and \cite{Surya} advocate simple methods of numerical Laplace inversion. Numerical computation of scale functions has already proved to be of practical value in, for example, the work of \cite{EM2008} and \cite{HR2002}.
\end{remark}

\section{Proof of Theorem \ref{rongsitu} in a subclass $\mathcal{S}^{(\infty)}\subseteq\mathcal{S}$}\label{piecewise}

In this section, our objective is to define a subclass $\mathcal{S}^{(\infty)}$ of $\mathcal{S}$ for which Theorem \ref{rongsitu} holds. To this end, 
by taking advantage of the fact that when $X$ has bounded variation, $0$ is irregular for $(-\infty,0)$, let us construct a pathwise solution to (\ref{SDE}) for $X$ having bounded variation and satisfying (H) (which will shortly turn out to be the unique solution within that class).
Define the   times $T_n$ and $S_n$ recursively as follows. We set $S_0=0$ and for $n=1,2,\ldots$
\begin{equation*}
\begin{split}
& T_n=\inf\{t>S_{n-1}:X_t-\delta\sum_{i=1}^{n-1}(S_i-T_i)\geq b\}, \\
& S_n=\inf\{t>T_n: X_t-\delta\sum_{i=1}^{n-1}(S_i-T_i) -\delta (t-T_n)<b\}.
\end{split}
\end{equation*}
Since $0$ is irregular for $(-\infty,0)$, the difference between two consecutive  times is strictly positive (except possibly for $S_0$ and $T_1$).
Now  we construct  a solution to (\ref{SDE}), $U=\{U_t : t\geq 0\}$, as follows. The process is issued from  $X_0=x$ and
\begin{equation*}
U_t =
\begin{cases}
X_t-\delta\sum_{i=1}^n(S_i-T_i) & \text{for $t\in[S_n,T_{n+1})$ and $n=0,1,2,\ldots$}\\
X_t-\delta\sum_{i=1}^{n-1}(S_i-T_i) -\delta (t-T_n) & \text{for $t\in[T_n,S_{n})$ and $n=1,2,\ldots$}
\end{cases}
\end{equation*}
Note that in particular the  times $T_n$ and $S_n$ for $n=1,2,\ldots$ can then be identified as
\begin{equation*}
T_n=\inf\{t>S_{n-1}: U_t \geq b\}, \quad 
S_n=\inf\{t>T_n: U_t <b\}
\end{equation*}
and moreover
\[
U_t = X_t - \delta\int_0^t \mathbf{1}_{\{U_s >b\}}{\rm d}s.
\]



The next Lemma is the first step in showing that any solution to (\ref{SDE}) which is not driven by a spectrally negative L\'evy process of bounded variation can be shown to exist uniquely as the result of of strong approximation by solutions to (\ref{SDE}) driven by a sequence of bounded variation processes respecting the condition (H).
In order to state it we shall introduce some notation.
\begin{definition}\rm
 It is known (cf. p.210 of \cite{Bert1996} for example) that for any spectrally negative L\'evy process with unbounded variation paths, $X$, there exists a sequence of bounded variation spectrally negative L\'evy processes, $X^{(n)}$, such that for each $t>0$,
\[
 \lim_{n\uparrow\infty}\sup_{s\in[0,t]}|X_s^{(n)} - X_s|=0
\]
almost surely and moreover, when $X^{(n)}$ is written in the form (\ref{sn}) the drift coefficient tends to infinity as $n\uparrow\infty$. The latter fact implies that for all $n$ sufficiently large, the sequence $X^{(n)}$ will automatically fulfil the condition (H). Such a sequence, $X^{(n)}$ will be referred to as {\it strongly approximating} for $X$. Rather trivially we may also talk of a strongly approximating sequence for processes of bounded variation respecting (H).
\end{definition}

\begin{lemma}\label{strongcgce} 
Suppose that $X$ is a spectrally negative L\'evy process satisfying (H) and that $X^{(n)}$ is a strongly approximating sequence.
Denote by $U^{(n)}$ the sequence of pathwise solutions associated with each $X^{(n)}$ which are constructed pathwise in the manner described above.
Then there exists a stochastic process $U^{(\infty)}=\{U^{(\infty)}_t : t\geq0\}$ such that for each fixed $t>0$, 
\[
\lim_{n\uparrow \infty}\sup_{s\in[0,t]}|U^{(n)}_s - U^{(\infty)}_s|=0
\]
almost surely.
\end{lemma}
{\bf Proof.} It suffices to give a proof for the case that $X$ has paths of unbounded variation.
Fix the constant $\eta>0$. Let $N\in\mathbb{N}$ be such that for all $n,m\geq N$, $\sup_{s\in[0,t]}|X^{(n)}(s)-X^{(m)}(s)|<\eta$ (note that in general $N$ is random). We will prove that for each fixed $t>0$
\begin{equation}
\label{localuniform}
\sup_{s\in[0,t]}|U^{(n)}_s-U_s^{(m)}|< 2 \eta
\end{equation}
from which we deduce that $\{U^{(n)}_s : s\in[0,t]\}$ is a Cauchy sequence in the Banach space consisiting of $D[0,t]$ equipped with the supremum norm where $D[0,t]$ is the space of cadlag mappings from $[0,t]$. Note that limit $U^{(\infty)}$ does not depend on $t$. Indeed, if $U^{(\infty,t_i)}$ for $i=1,2$ are the limits obtained over two different time horizons $0<t_1<t_2<\infty$ then a simple application of the triangle inequality shows that
\[
 \sup_{s\in[0,t_1]}|U^{(\infty,t_1)}_s-U_s^{(\infty, t_2)}|\leq  \lim_{n\uparrow\infty}\sup_{s\in[0,t_1]}|U^{(\infty, t_1)}_s-U_s^{(n)}| +\lim_{n\uparrow\infty}\sup_{s\in[0,t_1]}|U^{(\infty, t_2)}_s-U_s^{(n)}| =0
\]
almost surely.

Returning to the proof of (\ref{localuniform}),
define $\Delta^{(n,m)}U_s=U^{(n)}_s-U^{(m)}_s$ and $\Delta^{(n,m)} X_s=X^{(n)}_s-X^{(m)}_s$. Moreover, set
\begin{equation}
\label{drift_term}
A^{(n,m)}_s:=\Delta^{(n,m)} U_s-\Delta^{(n,m)} X_s=\delta\int_0^s \left( \mathbf{1}_{\{U^{(m)}_v>b,U^{(n)}_v\leq b\}} - \mathbf{1}_{\{U^{(m)}_v\leq b,U^{(n)}_v > b\}} \right) \mathrm{d}v.
\end{equation}
We shall proceed now to show that, almost surely  $\sup_{s\in[0,t]}|A^{(n,m)}_s|\leq\eta$ from which (\ref{localuniform}) follows.

Suppose the latter claim is not true. Then since $A^{(n,m)}$ is continuous
and $A^{(n,m)}_0=0$ there exists $0<s< t$  such that either (i) $A^{(n,m)}_s=\eta$ and for all $\epsilon>0$ sufficiently small there exists $r\in(s,s+\epsilon)$ such that $A^{(n,m)}_{r}>\eta$ or (ii) $A^{(n,m)}_s=-\eta$ and and for all $\epsilon>0$ sufficiently small there exists $r\in(s,s+\epsilon)$ such that $A^{(n,m)}_{r}<-\eta$.

In case (i) it follows that  $\Delta^{(n,m)} U_s>0$ since $\Delta^{(n,m)} X_s\in(-\eta,\eta)$ and thus by right-continuity there exists $\epsilon>0$ such that $\Delta^{(n,m)} U_r>0$ for all $r\in[s,s+\epsilon)$. Hence considering the integrand in \eqref{drift_term}, the first indicator is necessarily zero when $v\in [s,s+\epsilon)$. It follows that $A_{r}\leq \eta$ for all $r\in [s,s+\epsilon)$ which forms a contradiction.  
A similar argument by contradition excludes case (ii). 
\QED

We may now introduce the class $\mathcal{S}^{(\infty)}\subseteq\mathcal{S}$ for which we will be able to prove that the statement of Theorem \ref{rongsitu} holds. 

\begin{definition}
The class $\mathcal{S}^{(\infty)}=\mathcal{S}^{(\infty)}(x)$  consists of all processes $X\in\mathcal{S}$ (issued from $x$) such that for the associated process $U^{(\infty)}$ it holds that $\mathbb{P}_x(U^{(\infty)}_t=b)=0$ for Lebesgue almost every $t\geq 0$.
\end{definition}

\begin{remark}\rm\label{includesBV}
Note in particular that $\mathcal{S}^{(\infty)}$ contains all solutions to (\ref{SDE}) for which $X$ is of bounded variation satisfying (H).
\end{remark}

\begin{proposition}\label{almostS}
 When $X\in\mathcal{S}^{(\infty)}$, the process $U^{(\infty)}$ is the unique strong solution of \eqref{SDE} and consequently Theorem \ref{rongsitu} holds when the class $\mathcal{S}$ is replaced by $\mathcal{S}^{(\infty)}$.
\end{proposition}

{\bf Proof.} 
The fact that $U^{(\infty)}$ is a strong solution to (\ref{SDE}) is immediate as soon as it is clear that for each fixed $t>0$
\[
 \lim_{n\uparrow\infty}\int_0^t \mathbf{1}_{\{ U^{(n)} _s>b\}}\D s = \int_0^t \mathbf{1}_{\{ U^{(\infty)}_s >b\}}\D s
\]
almost surely. However this is an immediate consequence of Lemma \ref{strongcgce} and the assumption that $X\in \mathcal{S}^{(\infty)}$.

For pathwise uniqueness of this solution we use an argument which is based on ideas found in Example 2.4 on p286 of \cite{KS}. Suppose that $U^{(1)}$ and $U^{(2)}$ are two strong solutions to (\ref{SDE}) then writing $$\Delta_t = U^{(1)}_t- U^{(2)}_2 = -\delta\int_0^t (\mathbf{1}_{\{U^{(1)}_s >b \}  }    - \mathbf{1}_{\{U^{(2)}_s >b \}})\D s$$ it follows from classical calculus that 
\[
\Delta^2_t  = -2\delta\int_0^t \Delta_s (\mathbf{1}_{\{U^{(1)}_s >b \}  }    - \mathbf{1}_{\{U^{(2)}_s >b \}})\D s.
\]
Now note that thanks to the fact that $\mathbf{1}_{\{x>b\}}$ is an increasing function, it follows from the above representation that, for all $t\geq 0$, $\Delta^2_t\leq 0$ and hence $\Delta_t =0$ almost surely. This concludes the proof of existence and uniqueness amongst the class of strong solutions.  
\QED

\section{A Key analytical identity}

The main goal of this section is to establish a key analytical identity which will play an important role throughout the remainder of the paper.

\begin{theorem}\label{keyidentity} Suppose $X$ is a spectrally negative L\'evy process that has paths of  bounded variation and let $0<\delta<c$, where $c=\gamma+\int_{(0,1)}x\Pi(\mathrm{d}x)$. Then for $v\geq u> m\geq0$
\begin{multline}
\label{equation}
\int_0^\infty\int_{(z,\infty)}W^{(q)}(z-\theta +m)\Pi(\mathrm{d}\theta)\left[\frac{\mathbb{W}^{(q)}(v-m-z)}{\mathbb{W}^{(q)}(v-m)} \mathbb{W}^{(q)}(u-m)-\mathbb{W}^{(q)}(u-m-z)\right]\mathrm{d}z \\
= -\frac{\mathbb{W}^{(q)}(u-m)}{\mathbb{W}^{(q)}(v-m)} \left( W^{(q)}(v) + \delta\int_m^v\mathbb{W}^{(q)}(v-z)W^{(q)\prime}(z)\mathrm{d}z \right) \\
+W^{(q)}(u) + \delta\int_m^u\mathbb{W}^{(q)}(u-z)W^{(q)\prime}(z)\mathrm{d}z.
\end{multline}
\end{theorem}

{\bf Proof}.  
We denote  $p(x,\delta) = \e_x(e^{-q\kappa^+_a}\mathbf{1}_{\{\kappa^+_a<\kappa^-_0\}})$.
Suppose that $x\leq b$. Then by conditioning on $U$ until it passes above $b$, we have
\begin{equation}
p(x,\delta) 
= \mathbb{E}_x\left(\mathrm{e}^{-q\tau_b^+}\mathbf{1}_{\{\tau_0^->\tau_b^+\}}\right)p(b, \delta)
=\frac{W^{(q)}(x)}{W^{(q)}(b)}p(b, \delta).
\label{twosidedmodified}
\end{equation}
where in the last equality we have used ({\ref{twosided}}) from the Appendix.
Let now $x\geq b$ and $x\leq a$. 
Recall the process $Y=\{Y_t: t\geq 0\}$ where $Y_t = X_t-\delta t$ and denote by ${\rm P}_x$  the law of the process  $Y$ when issued from $x$ (with ${\rm E}_x$ as the associated expectation operator).
Using respectively that $0$ is irregular for $(-\infty, 0)$ for $Y$, (\ref{twosided}), the Strong Markov Property, \eqref{twosidedmodified}
and \eqref{two-sided-resolvent}, we have 
\begin{eqnarray}
p(x,\delta) &= & 
\E_x\left(\mathrm{e}^{-q\tau_a^+}\mathbf{1}_{\{\tau_0^->\tau_a^+\}}\mathbf{1}_{\{\tau_b^->\tau_a^+\}}\right)
+\E_x\left(\mathrm{e}^{-q\tau_a^+}\mathbf{1}_{\{\tau_0^->\tau_a^+\}}\mathbf{1}_{\{\tau_b^-<\tau_a^+\}}\right)\notag\\
&= &   \frac{\mathbb{W}^{(q)}(x-b)}{\mathbb{W}^{(q)}(a-b)}+\E_x\left(\mathrm{e}^{-q\tau_b^-}\mathbf{1}_{\{\tau_b^-
<\tau_a^+\}}\mathbb{E}_{U_{\tau_b^-}}\left(\mathrm{e}^{-q\tau_a^+}\mathbf{1}_{\{\tau_0^->\tau_a^+\}}\right)\right)\notag\\
&= &
\frac{\mathbb{W}^{(q)}(x-b)}{\mathbb{W}^{(q)}(a-b)}+\frac{p(b, \delta)}{W^{(q)}(b)}\E_x\left(\mathrm{e}^{-q\tau_b^-}\mathbf{1}_{\{\tau_b^-<\tau_a^+\}}W^{(q)}(Y_{\tau_b^-})\right)\notag\\
&= & \frac{\mathbb{W}^{(q)}(x-b)}{\mathbb{W}^{(q)}(a-b)}
+\frac{p(b, \delta)}{W^{(q)}(b)}
\int_0^{a-b}\int_{(y,\infty)}W^{(q)}(b+y-\theta)\notag\\
&& \cdot\left[\frac{\mathbb{W}^{(q)}(x-b)\mathbb{W}^{(q)}(a-b-y)}{\mathbb{W}^{(q)}(a-b)}-\mathbb{W}^{(q)}(x-b-y)\right]\Pi(\mathrm{d}\theta)\mathrm{d}y. \label{twosidedabove}
\end{eqnarray}
By setting $x=b$ in \eqref{twosidedabove} we can now get an explicit expression for $p(b,\delta)$ using that $\mathbb{W}^{(q)}(0)=1/({c}-\delta)$
\begin{multline}
p(b, \delta)= W^{(q)}(b)\bigg\{({c}-\delta)\mathbb{W}^{(q)}(a-b)W^{(q)}(b)\\
\left. -\int_0^{a-b}\int_{(y,\infty)}W^{(q)}(b+y-\theta)\mathbb{W}^{(q)}(a-b-y)\Pi(\mathrm{d}\theta)\mathrm{d}y\right\}^{-1}.
\label{*}
\end{multline}
We now start with the second step which concerns simplifying the term involving the double integral in above expression.
Noting that for $\delta=0$  (the case that there is no refraction) we have by \eqref{twosided} for all $x\geq 0$
\begin{equation}
p(b, 0) = \mathbb{E}_b\left(\mathrm{e}^{-q\tau_a^+}\mathbf{1}_{\{\tau_0^->\tau_a^+\}}\right)=\frac{W^{(q)}(b)}{W^{(q)}(a)},
\label{**} 
\end{equation}
it follows from (\ref{*}) and (\ref{**}) that
\begin{eqnarray}
\label{Wqexpr}
\lefteqn{\int_0^{a-b}\int_{(y,\infty)}W^{(q)}(y-\theta+b)W^{(q)}(a-b-y)\Pi(\mathrm{d}\theta)\mathrm{d}y}&&\notag\\
&&\hspace{5cm}={c} W^{(q)}(b)W^{(q)}(a-b)-W^{(q)}(a).
\end{eqnarray}
As $a\geq b$ is taken arbitrarily, we set $a=x$ in the above identity and take Laplace transforms from $b$ to $\infty$ of both sides of the above expression. Denote by $\mathcal{L}_b$ the operator which satisfies $\mathcal{L}_bf[\lambda]:=\int_b^\infty\mathrm{e}^{-\lambda x}f(x)\mathrm{d}x$. Let $\lambda>\Phi(q)$. For the left hand side of (\ref{Wqexpr}) we get by using Fubini's Theorem
\begin{multline*}
\int_b^\infty\mathrm{e}^{-\lambda x}\int_0^\infty\int_{(y,\infty)}W^{(q)}(y-\theta+b)W^{(q)}(x-b-y)\mathrm{d}y\Pi(\mathrm{d}\theta)\mathrm{d}x\\
=\frac{\mathrm{e}^{-\lambda b}}{\psi(\lambda)-q}\int_0^\infty\int_{(y,\infty)}\mathrm{e}^{-\lambda y}W^{(q)}(y-\theta+b)\Pi(\mathrm{d}\theta)\mathrm{d}y.
\end{multline*}
For the right hand side of (\ref{Wqexpr}) we get
\begin{eqnarray*}
\lefteqn{\int_b^\infty\mathrm{e}^{-\lambda x}\left(W^{(q)}(x-b){c} W^{(q)}(b)-W^{(q)}(x)\right)\mathrm{d}x} &&\\
&&\hspace{4cm}=\frac{\mathrm{e}^{-\lambda b}}{\psi(\lambda)-q} {c} W^{(q)}(b)-
\int_b^\infty\mathrm{e}^{-\lambda x}W^{(q)}(x)\mathrm{d}x
\end{eqnarray*}
and so
\begin{equation}
\int_0^\infty\int_{(y,\infty)}\mathrm{e}^{-\lambda y}W^{(q)}(y-\theta+b)\Pi(\mathrm{d}\theta)\mathrm{d}y 
={c} W^{(q)}(b)-(\psi(\lambda)-q)\mathrm{e}^{\lambda b}  
\mathcal{L}_bW^{(q)}[\lambda]\label{Wqlapl}
\end{equation}
for $\lambda>\Phi(q)$.
Our objective is now to use (\ref{Wqlapl}) to show that 
for $q\geq0$, for $x\geq b$, we have
\begin{multline}
\label{WqYexpr}
\int_0^\infty\int_{(y,\infty)}W^{(q)}(b+y-\theta)\Pi(\mathrm{d}\theta)\mathbb{W}^{(q)}(x-b-y)\mathrm{d}y\\
=-W^{(q)}(x)+({c}-\delta)W^{(q)}(b)\mathbb{W}^{(q)}(x-b)
-\delta\int_b^x\mathbb{W}^{(q)}(x-y)W^{(q)\prime}(y)\mathrm{d}y.
\end{multline}
The latter identity then implies the statement of the theorem.

The equality in  (\ref{WqYexpr}) follows by taking Laplace transforms on both sides in $x$. To this end note that by (\ref{Wqlapl}) it follows that the Laplace transform of the left hand side equals (for $\lambda>\varphi(q))$
\begin{eqnarray}
&&\int_b^\infty\mathrm{e}^{-\lambda x}\int_0^\infty\int_{(y,\infty)}W^{(q)}(b+y-\theta)\mathbb{W}^{(q)}(x-b-y)\Pi(\mathrm{d}\theta)\mathrm{d}y\mathrm{d}x\notag\\
&&\hspace{2cm}=\frac{\mathrm{e}^{-\lambda b}}{\psi(\lambda) - \delta  \lambda-q}\left({c} W^{(q)}(b)-(\psi(\lambda)-q)\mathrm{e}^{\lambda b}
\mathcal{L}_bW^{(q)}\right).\label{R}
\end{eqnarray}
Since $\mathcal{L}_b\left(\int_b^xf(x-y)g(y)\mathrm{d}y\right)[\lambda]=(\mathcal{L}_0f)[\lambda](\mathcal{L}_bg)[\lambda]$ and  $\mathcal{L}_bW^{(q)\prime}[\lambda]=\lambda \mathcal{L}_bW^{(q)}[\lambda]-\mathrm{e}^{-\lambda b}W^{(q)}(b)$ (which follows by integration by parts) it follows that the Laplace transform of the right hand side of (\ref{WqYexpr}) is equal to the right hand side of (\ref{R}). Hence (\ref{WqYexpr}) holds for almost every $x\geq b$. Because both sides of 
(\ref{WqYexpr}) are continuous in $x$, we have that (\ref{WqYexpr}) holds for all $x\geq b$. 
\QED

\begin{remark}\rm
Close examination  of the proof of the last theorem shows that may easily obtain the identity in Theorem \ref{two-sided} (i) for the case that $X$ has paths of bounded variation. Indeed plugging (\ref{equation}) into (\ref{twosidedabove}) and (\ref{*}) and then using (\ref{*}) in (\ref{twosidedmodified}) and (\ref{twosidedabove}) gives the required identity.

Note however that although this method, may be used to obtain other identities in the case of bounded variation paths, it is not sufficient to reach the entire family of identities presented in this paper which explains why the forthcoming line of reasoning does not necessarily appeal directly to the observation above.
\end{remark}

\section{Some calculations for resolvents}

In this section, we shall always take  $X$ to be of bounded variation  satisfyng (H). Recall for this class of driving L\'evy processes, we know that (\ref{SDE}) has a unique strong solution by Proposition \ref{almostS} which has been described piecewise at the beginning of Section \ref{piecewise}.

Define for $q>0$ and Borel $B\in[0,\infty)$,
\[
V^{(q)}(x,B)=\int_0^\infty \mathrm{e}^{-qt}\mathbb{P}_x(U_t\in B, \, \overline{U}_t\leq a, \, \underline{U}_t\geq 0)\mathrm{d}t = \int_0^\infty \mathbb{P}_x(U_t\in B, \, t<\kappa^-_0\wedge\kappa^+_a)  \mathrm{d}t.
\]
The identity in Theorem \ref{keyidentity} will be instrumental in establishing the following result.

\begin{proposition}\label{pseudo}When $X$ is of bounded variation satisfying (H)
the conclusion of Theorem \ref{dividends} (i) holds.
\end{proposition}

{\bf Proof.} Recall that the process $Y=\{Y_t: t\geq 0\}$ is given by $Y_t = X_t-\delta t$ and its law is  denote by ${\rm P}_x$  when issued from $x$.

We have for $x\leq b$  by the Strong Markov Property, \eqref{twosided} and \eqref{SNtwosidedresolve} 
\begin{eqnarray}
V^{(q)}(x,B)&= & \mathbb{E}_x\left(\int_0^{\tau_b^+}\mathrm{e}^{-qt}\mathbf{1}_{\{U_t\in B,t<\tau^+_a\wedge \tau^-_0\}}\mathrm{d}t \right)
+\mathbb{E}_x\left(\int_{\tau_b^+}^\infty\mathrm{e}^{-qt}\mathbf{1}_{\{U_t\in B,t<\tau^+_a\wedge \tau^-_0,\tau_b^+<\tau_0^-\}}\mathrm{d}t \right)\notag\\
&= & \mathbb{E}_x\left(\int_0^{\tau_b^+\wedge\tau_0^-}\mathrm{e}^{-qt}\mathbf{1}_{\{X_t\in B\}}\mathrm{d}t \right)
+\mathbb{E}_x\left( \mathrm{e}^{-q\tau_b^+}\mathbf{1}_{\{\tau_b^+<\tau_0^-\}}\right) V^{(q)}(b,B) \notag\\
&= & \int_B\left(\frac{W^{(q)}(b-y)}{W^{(q)}(b)}W^{(q)}(x)-W^{(q)}(x-y)\right)\mathrm{d}y+\frac{W^{(q)}(x)}{W^{(q)}(b)}V^{(q)}(b,B).
\label{x<b}
\end{eqnarray}
Moreover, for $b\leq x\leq a$ we have
\begin{eqnarray*}
\lefteqn{V^{(q)}(x,B)} &&\\
 &&= \E_x\left(\int_0^{\tau_{b}^-\wedge\tau_a^+}\mathrm{e}^{-qt}\mathbf{1}_{\{Y_t\in B\cap[b,a]\}}\mathrm{d}t \right)
+\E_x\left(\mathbf{1}_{\{\tau_{b}^-<\tau_a^+\}}\int_{\tau_{b}^-}^{\tau^+_a\wedge \tau^-_0}\mathrm{e}^{-qt}\mathbf{1}_{\{U_t\in B\}}\mathrm{d}t \right)\\
&&= \int_0^{\infty}\mathrm{e}^{-qt}\p_x\left(Y_t\in B\cap[b,a],t<\tau_{b}^-\wedge\tau_a^+\right)\mathrm{d}t
+\E_x\left( \mathbf{1}_{\{\tau_{b}^-<\tau_a^+\}}\mathrm{e}^{-q\tau_{b}^-}V^{(q)}(Y_{\tau_{b}^-},B)\right) \\
&&= \int_{B\cap[b,a]}\left(\frac{\mathbb{W}^{(q)}(a-z)}{\mathbb{W}^{(q)}(a-b)} \mathbb{W}^{(q)}(x-b)-\mathbb{W}^{(q)}(x-z)\right)\mathrm{d}z\\
&& + \int_0^\infty\int_{(z,\infty)} \Bigg\{ \int_B \left[ \frac{W^{(q)}(b-y)}{W^{(q)}(b)}W^{(q)}(z-\theta +b)-W^{(q)}(z-\theta+b-y) \right] \mathrm{d}y \\
&& + \frac{V^{(q)}(b,B)}{W^{(q)}(b)}W^{(q)}(z-\theta +b) \Bigg\} \left[ \frac{\mathbb{W}^{(q)}(a-b-z)}{\mathbb{W}^{(q)}(a-b)} \mathbb{W}^{(q)}(x-b)-\mathbb{W}^{(q)}(x-b-z)\right]\Pi(\mathrm{d}\theta)\mathrm{d}z
\end{eqnarray*}
where in the second equality we have used the Strong Markov Property and in the third equality \eqref{x<b}, (\ref{SNtwosidedresolve}) and (\ref{two-sided-resolvent}).
Next we shall apply the identity proved in Theorem \ref{keyidentity}  twice in order to simplify the expression for $V^{(q)}(x,B)$, $a\geq x\geq b$. We use it once by setting $m=b$, $u=x$, $v=a$ and once by setting $m=b-y$ and $u=x-y$, $v=a-y$ for $y\in[0,b]$. One obtains
\begin{equation}
\begin{split}
\label{x>b}
V^{(q)}(x,B) = & \int_{B\cap[b,a]}\left(\frac{\mathbb{W}^{(q)}(a-z)}{\mathbb{W}^{(q)}(a-b)} \mathbb{W}^{(q)}(x-b)-\mathbb{W}^{(q)}(x-z)\right)\mathrm{d}z  \\
& + \int_{B\cap[0,b)} \Bigg\{ \frac{W^{(q)}(b-y)}{W^{(q)}(b)}\bigg( -\frac{\mathbb{W}^{(q)}(x-b)}{\mathbb{W}^{(q)}(a-b)} \left( W^{(q)}(a) + \delta\int_{b}^a\mathbb{W}^{(q)}(a-z)W^{(q)\prime}(z)\mathrm{d}z \right)  \\
& +W^{(q)}(x) + \delta\int_{b}^x\mathbb{W}^{(q)}(x-z)W^{(q)\prime}(z)\mathrm{d}z \bigg)  \\
& - \bigg( -\frac{\mathbb{W}^{(q)}(x-b)}{\mathbb{W}^{(q)}(a-b)} \left( W^{(q)}(a-y) + \delta\int_{b-y}^{a-y}\mathbb{W}^{(q)}(a-y-z)W^{(q)\prime}(z)\mathrm{d}z \right)  \\
& +W^{(q)}(x-y) + \delta\int_{b-y}^{x-y}\mathbb{W}^{(q)}(x-y-z)W^{(q)\prime}(z)\mathrm{d}z \bigg)  \Bigg\} \mathrm{d}y  \\
& + \frac{V^{(q)}(b,B)}{W^{(q)}(b)} \bigg( -\frac{\mathbb{W}^{(q)}(x-b)}{\mathbb{W}^{(q)}(a-b)} \left( W^{(q)}(a) + \delta\int_{b}^a\mathbb{W}^{(q)}(a-z)W^{(q)\prime}(z)\mathrm{d}z \right)  \\
& +W^{(q)}(x) + \delta\int_{b}^x\mathbb{W}^{(q)}(x-z)W^{(q)\prime}(z)\mathrm{d}z \bigg).
\end{split}
\end{equation}
Setting $x=b$ in \eqref{x>b}, one then gets an expression for $V^{(q)}(b,B)$ in terms of itself. Solving this and then putting the resulting expression for $V^{(q)}(b,B)$ in \eqref{x<b} and \eqref{x>b}  leads to \eqref{2sidedres} which proves the proposition. 
\QED


Keeping with the setting that $X$ is a L\'evy process of bounded variation fulfilling (H), we may proceed to use the conclusion of the above proposition to establish an identity for the resolvent of $U$ without killing which we denote by 
\begin{equation}
 R^{(q)}(x, B)=\int_0^\infty e^{-q t}\mathbb{P}_x(U_t \in B)\D t
 \label{bounded by 1/q}
\end{equation}
for $q>0$ and Borel $B\in\mathbb{R}$.
\begin{corollary}
\label{corol_resolvent}
The conclusion of Theorem \ref{dividends} (iv) is valid when $X$ has paths of bounded variation satisfying (H).
\end{corollary}
{\bf Proof.}
By taking the expression given in Proposition \ref{pseudo} and letting $a\uparrow\infty$ one gets by the Monotone Convergence Theorem the expression for the one sided exit below resolvent given in Theorem \ref{dividends} (ii) in case $X$ is of bounded variation. It should be noted that here one uses the relation (cf. Chapter 8 of \cite{K2006})  $W^{(q)}(x)=\mathrm{e}^{\Phi(q)x}W_{\Phi(q)}(x)$ (and similarly $\mathbb{W}^{(q)}(x)=\mathrm{e}^{\varphi(q)}\mathbb{W}_{\varphi(q)}(x)$), where $\Phi(q)$ is the right inverse of the Laplace exponent  of $X$ and $W_{\Phi(q)}$  plays the role of the ($q=0$) scale function for the spectrally negative L\'evy process with Laplace exponent  $\psi(\theta+\Phi(q))-q$. Moreover one should use the known fact that $W_{\Phi(q)}(\infty)<\infty$ when $q>0$.

In the same spirit, replacing $b$ by $b+\theta$, $x$ by $x+\theta$, $B$ by $B+\theta$ and then letting $\theta\uparrow\infty$ in the expression for the one sided exit below resolvent obtained above, one may recover the expression for $ R^{(q)}(x, B)$ given in Theorem \ref{dividends} (iv). 
Here one uses L'H\^opital's rule and the known fact that the (left-) derivative of $W_{\Phi(q)}$ is bounded on intervals of the form $(x_0, \infty)$ where $x_0>0$ and tends to zero at infinity.

\QED

We close this section with a result which says that if $X^{(n)}$ strongly approximates $X$ and the latter has unbounded variation, then  by taking $n\uparrow\infty$, even though we are not yet able to necessarily identify the limiting process $U^{(\infty)}$ as the solution to (\ref{SDE}), the limit of the associated resolvents to $U^{(n)}$ say $R^{(q)}_n$ still exists as $n\uparrow\infty$ and it is absolutely continuous with density which is equal to the limiting density of $R^{(q)}_n$.

\begin{lemma}\label{throughintegral} Suppose that $X$ has paths of unbounded variation with strongly approximating sequence $X^{(n)}$. For $ x\in\mathbb{R}$  and bounded interval $B$ we have
\[
\lim_{n\uparrow\infty}R^{(q)}_n(x,B) = \int_B \lim_{n\uparrow \infty} r_n^{(q)}(x,y)\D y,
\]
where $r_n^{(q)}(x,y)$ is the density of $R^{(q)}_n(x,\mathrm{d}y)$.
In particular 
 $\lim_{n\uparrow\infty} r_n^{(q)}(x,y)$ is equal to the density in the right hand side of  (\ref{2sidedres}).
\label{limintintlim}
\end{lemma}
\textbf{Proof}
The proof is a direct consequence of the Dominated Convergence Theorem and the fact that, by the Continuity Theorem for Laplace transforms, both $W^{(q)}$, $\mathbb{W}^{(q)}$ and $W^{(q)\prime}$ are continuous with respect to the L\'evy triplet of the underlying L\'evy process. 
\QED

\section{Proof of Theorem \ref{rongsitu}}

Our objective in this section is to use the resolvents of the previous section to prove the following result.
\begin{lemma}\label{S0=S}
It holds that $\mathcal{S}^{(\infty)}$ contains all spectrally negative L\'evy processes of unbounded variation and hence by Remarks \ref{strongsigma} and \ref{includesBV}  and Proposition \ref{almostS} it follows that Theorem \ref{rongsitu} holds. 
\end{lemma}

{\bf Proof.}
It suffices to show that for all driving L\'evy processes $X$ with paths of unbounded variation we have that, when $x$ is fixed, $\mathbb{P}_x(U^{(\infty)}_t=b)=0$ for Lebesgue almost every $t\geq 0$.
In fact we shall prove something slightly more general (for future convenience). 

Let $X$ be strongly approximated by the sequence $X^{(n)}$. Note that for each $t,\eta>0$ and $a\in\mathbb{R}$, thanks to Lemma \ref{strongcgce}.
\[
\{U^{(\infty)}_t =a\}\subseteq\liminf_{n\uparrow \infty }\{U^{(n)}_t \in (a-\eta, a+\eta)\}:=\{U^{(n)}_t \in (a-\eta, a+\eta) \text{ eventually as } n\uparrow \infty \}.
\]
Standard measure theory (cf. Exercise 3.1.12 of \cite{Strook}) now gives  us  for each $\eta>0$
\[
\mathbb{P}_x(U^{(\infty)}_t =a)\leq \liminf_{n\uparrow \infty} \mathbb{P}_x(U^{(n)}_t \in (a-\eta, a+\eta)). 
\]
Now applying Fatou's Lemma followed by the conclusion of Lemma \ref{limintintlim} we have for 
 $q>0$,
 \begin{eqnarray*}
 \int_0^\infty e^{-qt}\mathbb{P}_x(U^{(\infty)}_t =a)\D t&\leq &
 \liminf_{n\uparrow \infty}  \int_0^\infty e^{-q t}\mathbb{P}_x(U^{(n)}_t \in (a-\eta, a+\eta))\D t\\
& =&\liminf_{n\uparrow \infty}R^{(q)}_n(x,(a-\eta, a+\eta))\\
&=& \int_{a-\eta}^{a+\eta} r^{(q)}(x,y)\D y
 \end{eqnarray*}
 where $ r^{(q)}(x,y)$ is the density on the right hand side of   (\ref{2sidedres}). Note that, uniformly in $\eta$, the integral on the right hand side above is bounded by $1/q$ thanks to  Lemma \ref{throughintegral} and the fact that for all $n$, $R^{(q)}_n(x,\mathbb{R})\leq 1/q$ on account of (\ref{bounded by 1/q}).
 Since the quantity $\eta$ is arbitrary the required statement that $\mathbb{P}_x(U^{(\infty)}_t=a)=0$ for Lebesgue almost every $t>0$ follows. \QED

Before concluding this section, it is worth registering the following corollary for the next section which follows directly from the conclusion and proof above.

\begin{corollary}\label{nextsectionlemma}
For all $X\in\mathcal{S}$, we have for each given $x,a\in\mathbb{R}$ that the unique strong solution $U$ to (\ref{SDE}) satisfies $\mathbb{P}_x(U_t = a)=0$ for Lebesgue almost every $t\geq 0$.
\end{corollary}

\section{Proof of Theorem \ref{dividends}}

Firstly let us note that parts (ii), (iii) and (iv) follow from part (i) by taking limits much in the spirit of the proof of Corollary \ref{corol_resolvent}. As before such calculations are straightforward and hence, for the sake of brevity, are left to the reader.

To establish part (i) we have already seen that (\ref{2sidedres}) is true for case that $X$ has bounded variation and satisfies (H).
To deal with the case that $X$ has paths of unbounded variation we consider as usual a strongly approximating sequence $X^{(n)}$.
As before, we will write the left hand side of the identity in (\ref{2sidedres}) when the driving process is $X^{(n)}$ in the form
\[
V^{(q)}_n(x, B) =\int_0^\infty e^{-qt}\mathbb{P}_x(U^{(n)}_t \in B, \, \overline{U}^{(n)}_t \leq a, \, \underline{U}^{(n)}_t \geq 0) \D t
\]
 for $q>0$ and Borel $B\in \mathbb{R}$. 
 
 Recall from Lemma \ref{S0=S} that $U^{(\infty)}$ defined in Lemma \ref{strongcgce} is the unique solution to (\ref{SDE}). We shall henceforth refer to it as just $U$.  In the spirit of Lemma \ref{throughintegral} we may prove that for open intervals $B$,
\begin{equation}
\lim_{n\uparrow \infty} V_n^{(q)}(x, B)  = \int_B v^{(q)}(x,y)\D y
\label{referback}
\end{equation}
where $v^{(q)}(x, y)$ is the density which appears on the right hand side of the identity (\ref{2sidedres}). It is known (see for example Lemma 13.4.1  of \cite{whitt}) that 
\[
|\overline{U}^{(n)}_t - \overline{U}_t|\vee|\underline{U}^{(n)}_t - \underline{U}_t|\leq \sup_{s\in[0,t]}|U^{(n)}_s  - U_s|
\]
Thanks to Lemma \ref{strongcgce}, it follows that for each $t>0$, in the almost sure sense,
\[
\lim_{n\uparrow \infty}(U^{(n)}_t, \overline{U}^{(n)}_t,\underline{U}^{(n)}_t) = (U_t, \overline{U}_t, \underline{U}_t).
\]
This tells us that by the Dominated Convergence Theorem
\[
\lim_{n\uparrow \infty} V_n^{(q)}(x, B)=\int_0^\infty e^{-qt}\mathbb{P}_x(U_t\in B, \overline{U}_t \leq a, \underline{U}_t\geq 0)\D t
\]
{\it providing} the boundary of $\{U_t \in B, \, \overline{U}_t \leq a, \, \underline{U}_t \geq 0\}$  is not charged by $\mathbb{P}_x$. 
To rule the latter out it suffices to show that 
\begin{equation}
\mathbb{P}_x(U_{t} \in\partial B)=\mathbb{P}_x(\underline{U}_{t} =0)=\mathbb{P}_x(\overline{U}_{t}=a)=0.
\label{threethings}
\end{equation}
for Lebesgue almost every $t\geq 0$.

To this end, note that if $\kappa^{[0,\epsilon)}=\inf\{t>0: U_t \in [0,\epsilon)\}$ where $\epsilon>0$ then it is easy to see from (\ref{SDE}) that 
on $\{\kappa^{[0,\epsilon)}<\infty\}$
\[
U_{\kappa^{[0,\epsilon)} + s} \leq   U_{\kappa^{[0,\epsilon)}} + \widetilde{X}_s
\]
where $\widetilde{X}$ is a copy of $X$ which is independent of $\{U_s: s\leq \kappa^{[0,\epsilon)}\}$.
Let $g(x,t)=\mathbb{P}_x(\inf_{s\leq t}X_s\geq 0)$ and note that it is increasing in $x$ and decreasing in $t$. Moreover, by regularity of $(-\infty,0)$ for $X$ we have that for each fixed  $t>0$,  $g(x,t)=0$. 
It follows that for all $\epsilon>0$
\begin{eqnarray}
\lefteqn{\mathbb{P}_x(\underline{U}_{t} =0)}\notag\\
&&\leq \mathbb{E}_x[\mathbf{1}_{\{\kappa^{[0,\epsilon)}\leq t\}}\mathbb{P}_{x}  (U_{\kappa^{[0,\epsilon)}}+\inf_{s\leq t-\kappa^{[0,\epsilon)} }\widetilde{X}_{s}\geq 0 |\mathcal{F}_{\kappa^{[0,\epsilon)}})]\notag\\
&&\leq \mathbb{E}_x[\mathbf{1}_{\{\kappa^{[0,\epsilon)}\leq t\}} g(\epsilon, t- \kappa^{[0,\epsilon)} )].
\label{U=0}
\end{eqnarray}
By monotonicity there exist $\kappa^{\{0\}}:=\lim_{n\uparrow \infty} \kappa^{[0,\epsilon)}$
 and the event $\{\kappa^{\{0\}} =t\}$ implies almost surely that $U_{t-} = 0 =U_t$ where the last equality follows on account of the fact that $t$ is a jump time with probability zero. 
Hence by dominated convergence
\[
\mathbb{P}_x(\underline{U}_{t} =0)\leq 
\mathbb{E}_x[\mathbf{1}_{\{\kappa^{\{0\}}<t\}} 
\lim_{n\uparrow \infty}g(\epsilon,  t- \kappa^{\{0\}} )
]
+
\mathbb{P}_x(U_t =0).
\]
The preceding remarks concerning $g(x,t)$ and the 
the conclusion of Corollary \ref{nextsectionlemma}
now imply that
$\mathbb{P}_x(\underline{U}_{t} =0)=0$ for Lebesgue almost every $t\geq 0$. A similar argument can be employed to show that  $\mathbb{P}_x(\overline{U}_{t} =a)=0$ for Lebesgue almost every $t\geq 0$. It is also a simple consequence of Corollary \ref{nextsectionlemma} that $\mathbb{P}_x(U_{t} \in\partial B)=0$ for Lebesgue almost every $t\geq 0$.
 Thus (\ref{threethings}) is satisfied
and referring back to (\ref{referback}) we see that the proof is compete.
\QED

\section{Proof of Theorems \ref{two-sided} and \ref{ruin}}

The idea of  the proofs is to make use of the identities in parts (i)--(iv) of Theorem \ref{dividends}.
We give only the important ideas of the proof as the details of the computations are straightforward and so left to the reader, again, for the sake of brevity. In doing so, one will need to make use of the following identity for $q,a\geq0$
\begin{equation*}
\delta\int_0^a\mathbb{W}^{(q)}(a-y)W^{(q)}(y)\mathrm{d}y=\int_0^a\mathbb{W}^{(q)}(y)\mathrm{d}y-\int_0^aW^{(q)}(y)\mathrm{d}y,
\end{equation*}
which can be proved by showing that the Laplace transforms on both sides are equal.

One obtains the result in Theorem \ref{ruin} (ii)  by noting that 
\[
\mathbb{E}_x\left( \mathrm{e}^{-q\kappa_0^-}\mathbf{1}_{\{\kappa_0^-<\infty\}}\right)=1-\mathbb{P}_x(\underline{U}_{\mathbf{e}_q} \geq 0) =1 - q\int_0^\infty e^{-qt} \mathbb{P}_x(U_t\in\mathbb{R}, t<\kappa^-_0)\D t.
\]
For the proof of Theorem \ref{two-sided} (i), it suffices to note that for $q>0$, by applying the Strong Markov Property, one has that
\[
\mathbb{P}_x(\underline{U}_{\mathbf{e}_q} \geq 0, \overline{U}_{\mathbf{e}_q}>a) = \e_x(e^{-q\kappa^+_a}\mathbf{1}_{\{\kappa^+_a<\kappa^-_0\}}) \mathbb{P}_a(\underline{U}_{\mathbf{e}_q} \geq 0) .
\]
The first and the last probabilities above can be obtained directly from the potential measures given in Theorem \ref{dividends} 
since 
\begin{eqnarray*}
\mathbb{P}_x(\underline{U}_{\mathbf{e}_q} \geq 0, \overline{U}_{\mathbf{e}_q}>a) &=& \mathbb{P}_x(\underline{U}_{\mathbf{e}_q} \geq 0)- \mathbb{P}_x(\underline{U}_{\mathbf{e}_q} \geq 0, \overline{U}_{\mathbf{e}_q}\leq a)\\
&=&q \int_0^\infty e^{-qt} \mathbb{P}_x(U_t\in\mathbb{R}, t<\tau^-_0) - q \int_0^\infty e^{-qt} \mathbb{P}_x(U_t\in[0,a], t<\tau^+_a \wedge \tau^-_0)\D t
\end{eqnarray*}
and 
\[
\mathbb{P}_a(\underline{U}_{\mathbf{e}_q} \geq 0)  = q \int_0^\infty e^{-qt} \mathbb{P}_a(U_t\in[0,\infty), t<\tau^-_0)\D t.
\]

By using the Strong Markov Property for (\ref{SDE}) at the specific stopping time $\kappa^+_a$ and the fact that $U_{\kappa^+_a} =a$ on $\{\kappa^+_a<\infty\}$ we now have that 
\begin{equation*}
\begin{split}
\mathbb{E}_x\left( \mathrm{e}^{-q\kappa_0^-}\mathbf{1}_{\{\kappa_0^-<\kappa_a^+\}} \right) = & \mathbb{E}_x\left( \mathrm{e}^{-q\kappa_0^-}\mathbf{1}_{\{\kappa_0^-<\infty\}} \right) - \mathbb{E}_x\left( \mathrm{e}^{-q\kappa_0^-}\mathbf{1}_{\{\kappa_a^+<\kappa_0^-\}} \right) \\
= & \mathbb{E}_x \left( \mathrm{e}^{-q\kappa_0^-}\mathbf{1}_{\{\kappa_0^-<\infty\}} \right) - \mathbb{E}_x \left( \mathrm{e}^{-q\kappa_a^+}\mathbf{1}_{\{\kappa_a^+<\kappa_0^-\}} \right) \mathbb{E}_a\left( \mathrm{e}^{-q\kappa_0^-}\mathbf{1}_{\{\kappa_0^-<\infty\}} \right),
\end{split}
\end{equation*}
for $0\leq x,b\leq a$. This gives the required identity in Theorem \ref{two-sided} (ii).

For part (i) of Theorem \ref{ruin} one notes that
\[
\mathbb{E}_x\left( \mathrm{e}^{-q\kappa_a^+}\mathbf{1}_{\{\kappa_a^+<\infty\}}\right) = 1-\mathbb{P}_x(\overline{U}_{\mathbf{e}_q} \leq a) = 1- q\int_0^\infty e^{-qt} \mathbb{P}_x(U_t\in(-\infty, a], t<\kappa^+_a)\D t.
\]
However, it seems difficult to derive the required expression in this way. In place of this method one may obtain the result by using the expression for  $\e_x(e^{-q\kappa^+_a}\mathbf{1}_{\{\kappa^+_a<\kappa^-_0\}})$, namely first replace $x$ by $x+\theta$, $a$ by $a+\theta$ and $b$ by $b+\theta$ and then let $\theta\uparrow\infty$.

\QED

\section{Proof of Theorem \ref{creep}}
It is a well established fact (cf. Chapter VI of \cite{Bert1996}) that a spectrally negative L\'evy process creeps downward if and only if it has a Gaussian component. For this reason it is obvious that the probability that $U$ creeps downward is zero as soon as $X$ has no Gaussian component. We therefore restrict ourselves to the case that $X$ has no Gaussian component.

Suppose that for $x,b\geq a$, $w^{(q)}(x,y,a,b)$ is the resolvent density for $U$ with killing on exiting the interval  $[a,\infty)$.  Note that by spatial homogeneity $w^{(q)}(x,y,a,b)=w^{(q)}(x-a,y-a,0,b-a)$ and therefore an expression for this density is already given in Theorem \ref{dividends}. Since $U$ is the sum of a continuous process and a L\'evy process, it is quasi-left continuous and hence we can use 
 Proposition 1(i) in \cite{pistorius_potential} to deduce 
\begin{equation*}
\begin{split}
\mathbb{E}_x\left( \mathrm{e}^{-q\kappa_0^-}\mathbf{1}_{\{ U_{\kappa_0^-}=0 \}}   \right) = & \lim_{\epsilon\downarrow 0}\mathbb{E}_x \left( \mathrm{e}^{-q\kappa^{\{0\}}}\mathbf{1}_{\{  \kappa^{\{0\}}<\kappa_{-\epsilon}^-  \}}  \right) = \lim_{\epsilon\downarrow 0} \frac{w^{(q)}(x,0,-\epsilon,b)} {w^{(q)}(0,0,-\epsilon,b)} \\
= & \lim_{\epsilon\downarrow 0} \frac{w^{(q)}(x+\epsilon,\epsilon,0,b+\epsilon)} {w^{(q)}(\epsilon,\epsilon,0,b+\epsilon)},
\end{split}
\end{equation*}
where $\kappa^{\{0\}}=\inf\{t>0:U_t=0\}$.
The limit can then be computed by using l'H\^optital's rule, the Dominated Convergence Theorem and the fact that $W^{(q)\prime}(0)=2/\sigma^2$ when $\sigma>0$.
\QED

\section{Applications in ruin theory}\label{applications}

As alluded to in the introduction, modern perspectives on the theory of ruin has seen preference for working with spectrally negative L\'evy processes.
Indeed one may understand the third bracket in (\ref{generalCR})  as the part of a risk process corresponding to countably infinite number of arbitrarily small claims compensated by a deterministic positive drift (which may be infinite in the case that
$\int_{(0,1)} x\Pi(\D x)=\infty$) corresponding to the accumulation of premiums over an infinite number of contracts. Roughly speaking, the way in which claims occur is such that in any arbitrarily small period of time $\D t$, a claim of size $x$ is made independently with probability $\Pi(\D x)\D t + o(\D t)$. The insurance company thus counterbalances
such claims by ensuring that it collects premiums in such a way that in any $\D t$, $x\Pi(\D x)\D t$ of its income is devoted to the compensation of claims of size $x$.
The second bracket in (\ref{generalCR}) we may understand as coming from large claims which occur occasionally and are compensated against by a steady income at rate $\gamma>0$ as in the Cram\'er-Lundberg model. Here `large' is taken to mean claims of size one or more.
Finally the first bracket in (\ref{generalCR}) may be seen as a stochastic perturbation of the system of claims and premium income.


As mentioned earlier, a quantity which is  of particular value is the probability of ruin. This is given precisely in the second half of Theorem \ref{ruin} (i).
Another quantity of interest mentioned in the introduction is the net present value of the dividends paid out until ruin.  Such a quantity is easily obtained from Theorem \ref{dividends} and it is equal to 
\begin{eqnarray}\label{Vdivs}
\mathbb{E}_x\left(\int_0^{\kappa^-_0} e^{-qt}\delta\mathbf{1}_{\{U_t >b\} } \D s \right)
&=&\frac{\delta}{q}\left(1-\mathbb{Z}^{(q)}(x-b)\right)\notag\\
&&+\frac{W^{(q)}(x)+\delta\mathbf{1}_{\{x\geq b\}} \int_b^x\mathbb{W}^{(q)}(x-y)W^{(q)\prime}(y)\mathrm{d}y}
{\varphi(q)\int_0^\infty\mathrm{e}^{-\varphi(q) y}W^{(q)\prime}(y+b)\mathrm{d}y}.
\end{eqnarray}
As $U$ is a semi-martinagle whose jumps are described by the same Poisson point process of jumps which describes the jumps of the driving L\'evy process, one may apply the compensation formula in a straightforward way together with the resolvent in part (ii) of Theorem \ref{dividends} to deduce the following expression for the joint law of the overshoot and undershoot at ruin (see for example the spirit of the discussion at the beginning of Section 8.4 of \cite{K2006}) in the case that 
$0<\delta<\mathbb{E}(X_1)$.

Let $A\subset(-\infty,0)$ and $B\subset[0,\infty)$ be Borel-sets and let $U_{\kappa^-_0-}=\lim_{t\uparrow\kappa^-_0}U_{t}$. 
For $x\in\mathbb{R}$ 
\begin{eqnarray*}
\lefteqn{\mathbb{P}_x(U_{\kappa^-_0}\in A,U_{\kappa^-_0-}\in B)} &&  \\
&&=
\int_B\Pi(y-A)\frac{1-\delta W(b-y)}{1-\delta W(b)}\mathrm{d}y\left(W(x)+\delta 
\mathbf{1}_{\{x\geq b\}}\int_b^x \mathbb{W}(x-z)W'(z)\mathrm{d}z\right) \\
&& - \int_{B\cap[0,b)}\Pi(y-A)\left( W(x-y)+\delta 
\mathbf{1}_{\{x\geq b\}}\int_b^x\mathbb{W}(x-z)W'(z-y)\mathrm{d}z\right)\mathrm{d}y \\
&& - \int_{B\cap[b,\infty)}\Pi(y-A)\mathbb{W}(x-y)\mathrm{d}y.
\end{eqnarray*}

As mentioned in the introduction, expressions for the expected discounted value of the dividends, the Laplace transform of the ruin probability and the joint law of the undershoot and overshoot have been established before for refracted L\'evy processes, but only for the case $\Pi(0,\infty)<\infty$. Moreover, the identities we have obtained here, aside from being more generally applicable, arguably appear in a simpler form, being expressed in terms of scale functions. 
For example, considering the expression for the value of the dividends, denoted by $V(x)$, given in \eqref{Vdivs}, we see that we can easily differentiate that expression with respect to $x$ (providing $W^{(q)},\mathbb{W}^{(q)}\in C^1(0,\infty)$). In that case, it follows that  there is smooth pasting, i.e. $\lim_{x\uparrow b}V'(x)=\lim_{x\downarrow b}V'(x)$, if and only if $X$ has paths of unbounded variation  or $b$ is chosen such that $\varphi(q)\int_0^\infty\mathrm{e}^{-\varphi(q) y}W^{(q)\prime}(y+b)\mathrm{d}y=W^{(q)\prime}(b)$. Having an expression for the derivative of $V$ is very important regarding a certain optimal control problem involving refracted L\'evy processes, see Gerber and Shiu \cite{GS2006} who solve this control problem for an extremely particular example of a refracted L\'evy process (the compound Poisson case with exponentially distributed jumps).
Besides in the Cram\'er-Lundberg model, the threshold strategy (and/or corresponding control problem) has also been considered in a Brownian motion setting, see e.g. \cite{AT1997, GS2006BM, JS1995}. Refracted L\'evy processes have also been recently studied in the context of queuing theory, see e.g. Bekker et al. \cite{bekkerboxmaresing} and references therein.  By comparison, the setting here operates at a greater degree of generality however.

\bigskip

\noindent We conclude this section with two concrete examples.


\subsection*{Example 1}
Suppose that we take $X$ to be a spectrally negative $\alpha$-stable process for $\alpha\in(1,2)$ with positive linear drift $c>\delta$. It is known that for such processes (cf. \cite{F1998}), 
\[
W(x) =\frac{1}{c}\left(1-E_{\alpha-1}(-c x^{\alpha-1})\right)
\]
where $E_{\alpha-1}(x) = \sum_{n\geq 0} x^n/\Gamma((\alpha-1)n+1)$ is the one-parameter Mittag-Leffer function with index $\alpha-1$. 
It follows that when $X$ is refracted with rate $\delta$, then the ruin probability is given by
\begin{equation*}
\begin{split}
 \mathbb{P}_x(\kappa^-_0<\infty)   
  = & 1 - \frac{c - \delta }{c - \delta + \delta E_{\alpha-1}(-c b^{\alpha-1})} \Big\{ 1-E_{\alpha-1}(-c x^{\alpha-1}) \\
& -\mathbf{1}_{\{x\geq b\}}  \delta(\alpha-1) \int_b^x [1 - E_{\alpha-1}(-(c-\delta) (x-y)^{\alpha-1})]    E'_{\alpha-1}(-c y^{\alpha-1}) y^{\alpha-2}\mathrm{d}y \Big\}.
\end{split}
\end{equation*}

\subsection*{Example 2}

Let $X$ be a spectrally  negative L\'evy  process of bounded variation with compound Poisson jumps
such that the L\'evy measure is given by
\begin{equation*}
\Pi(\mathrm{d}x)=\lambda\sum_{k=1}^n A_k\mathrm{e}^{-\alpha_k x}\mathrm{d}x, \quad \lambda,A_k,\alpha_k>0,\sum_{k=1}^nA_k=1
\end{equation*}
and when written in the form (\ref{sn}) the drift coefficient is taken to be $c$ such that $\mathbb{E}(X_1)>0$. 
This corresponds to the case of a Cram\'er-Lundberg process with premium rate $c$ and claims which are hyper-exponentially distributed. Moreover we assume that $q>0$ and that $0<\delta<c$. Then the Laplace exponent of $X$ is well defined and given by
\begin{equation*}
\log \mathbb{E}\left( \mathrm{e}^{\theta X_1} \right) = c\theta - \lambda + \lambda\sum_{k=1}^n A_k\frac{\alpha_k}{\alpha_k+\theta} \quad \text{for $\theta>\min\{\alpha_1,\ldots,\alpha_n\}$}.
\end{equation*}
Denote (with slight abuse of notation) by $\psi(\theta)$ the right hand side of above equation and note that this expression is well defined for all $\theta\in\mathbb{R}\backslash\{-\alpha_1,\ldots,-\alpha_n\}$.
By using the partial fraction method, we can then write for all $\theta\in\mathbb{R}\backslash\{-\alpha_1,\ldots,-\alpha_n\}$, 
\begin{equation*}
\begin{split}
\frac{1}{\psi(\theta)-q} = & \frac{1}{c\theta - \lambda + \lambda\sum_{k=1}^n A_k\frac{\alpha_k}{\alpha_k+\theta} - q } \cdot\frac{\prod_{k=1}^n(\alpha_k+\theta) }
{\prod_{k=1}^n(\alpha_k+\theta)} \\
= & \frac{\prod_{k=1}^n(\alpha_k+\theta)}{c\prod_{i=0}^n(\theta-\theta_i)}=\sum_{i=0}^n\frac{D_i}{\theta-\theta_i}.
\end{split}
\end{equation*}
Here $\{\theta_i :  i =0,1,...,n \}$ are the roots of $\psi(\theta)-q$, with $\theta_0=\Phi(q)>0$ and the other roots being strictly negative. Further $\{D_i : i =0,1,...,n\}$ are given by $D_i=1/\psi'(\theta_i)$. It follows that the scale function of $X$ is given by
\begin{equation*}
W^{(q)}(x)=\sum_{i=0}^n D_i\mathrm{e}^{\theta_i x }, \quad x\geq0.
\end{equation*}
Similarly, the scale function of the process $\{X_t-\delta t: t\geq0\}$ is given by
\begin{equation*}
\mathbb{W}^{(q)}(x)=\sum_{j=0}^n \tilde{D}_j\mathrm{e}^{\tilde{\theta}_j x}, \quad x\geq0,
\end{equation*}
where $\{\tilde{\theta}_j: j=0,1,...,n \}$ are the roots of $\psi(\theta)-\delta\theta-q$ with $\tilde{\theta}_0=\varphi(q)>0$ and $\tilde{D}_j=1/\psi'(\tilde{\theta}_j)$. We now want to give an explicit expression for the value of the dividends, denoted by $V$, for which the generic formula was given in (\ref{Vdivs}) above.

We can write
\begin{equation*}
\begin{split}
\int_{b}^{x}\mathbb{W}^{(q)}(x-z)W^{(q)\prime}(z)\mathrm{d}z = & \sum_{j=0}^n\sum_{i=0}^n \frac{\tilde{D}_j }{\theta_i-\tilde{\theta}_j} D_i \theta_i
\left( \mathrm{e}^{\theta_i x}-\mathrm{e}^{\theta_i b}\mathrm{e}^{\tilde{\theta}_j(x-b)}   \right) \\
= & -\frac{1}{\delta}W^{(q)}(x) - \sum_{j=0}^n\sum_{i=0}^n \frac{\tilde{D}_j }{\theta_i-\tilde{\theta}_j} D_i \theta_i
 \mathrm{e}^{\theta_i b}\mathrm{e}^{\tilde{\theta}_j(x-b)},
\end{split}
\end{equation*}
where the second equality follows since
\begin{equation*}
\sum_{j=0}^n \frac{\tilde{D}_j }{\theta_i-\tilde{\theta}_j} = \frac{1}{\psi(\theta_i)-\delta\theta_i-q} = -\frac{1}{\delta\theta_i}.
\end{equation*}
Further we have
\begin{equation*}
\varphi(q)\int_0^\infty\mathrm{e}^{-\varphi(q)y}W^{(q)\prime}(y+b)\mathrm{d}y = \tilde{\theta}_0 \sum_{i=0}^n
 \frac{D_i \theta_i}{\tilde{\theta}_0 -\theta_i }\mathrm{e}^{\theta_i b}
\end{equation*}
and since $\sum_{j=0}^n \tilde{D}_j/\tilde{\theta}_j= -1/(\psi(0)-\delta\cdot0-q)$, we get
\begin{equation*}
\frac{\delta}{q} \left( 1-\mathbb{Z}^{(q)}(x-b)  \right) = -\delta\sum_{j=0}^n \frac{\tilde{D}_j}{\tilde{\theta}_j}
\left( \mathrm{e}^{\tilde{\theta}_j(x-b)} -1  \right) = \frac{\delta}{q} -\delta\sum_{j=0}^n \frac{\tilde{D}_j}{\tilde{\theta}_j}
  \mathrm{e}^{\tilde{\theta}_j(x-b)}.
\end{equation*}
Hence the value of the dividends $V$ is given for $x\leq b$ by
\begin{equation*}
V(x)=
\left( \tilde{\theta}_0 \sum_{i=0}^n \frac{D_i \theta_i}{\tilde{\theta}_0 -\theta_i }\mathrm{e}^{\theta_i b} \right)^{-1}\cdot \sum_{i=0}^n D_i\mathrm{e}^{\theta_i x }    
\end{equation*}
and for $x\geq b$ by
\begin{equation*}
V(x)=  \frac{\delta}{q} +  \sum_{j=0}^n \left\{  \left( \tilde{\theta}_0 \sum_{i=0}^n \frac{D_i \theta_i}{\tilde{\theta}_0 -\theta_i }\mathrm{e}^{\theta_i b} \right)^{-1}  \sum_{i=0}^n \frac{\tilde{D}_j }{\tilde{\theta}_j-\theta_i} D_i \theta_i \mathrm{e}^{\theta_i b}  -  \frac{\tilde{D}_j}{\tilde{\theta}_j}  \right\} \delta\mathrm{e}^{\tilde{\theta}_j(x-b)}.
\end{equation*}
Note that the $j=0$ term between the curly brackets is zero. The above formulas for $V$ are an improvement upon the calculations made in Appendix A of Gerber \& Shiu \cite{GS2006}.

\section*{Appendix}


The theorem below is a collection of known fluctuation identities which have been used in the preceding text. See for example Chapter 8 of \cite{K2006} for proofs and the origin of these identities.

\begin{theorem}\label{appendixthrm}Recall that $X$ is a spectrally negative L\'evy process and let
\[
\tau^+_a = \inf\{t>0: X_t >a\}\text{ and }\tau^-_0=\inf\{t>0 : X_t <0\}.
\]
\begin{itemize}
\item[(i)] For $q\geq0$ and $x\leq a$
\begin{equation}
\label{twosided}
\mathbb{E}_x\left(\mathrm{e}^{-q\tau_a^+}\mathbf{1}_{\{\tau_0^->\tau_a^+\}}\right)=\frac{W^{(q)}(x)}{W^{(q)}(a)}.
\end{equation}
\item[(ii)] For any $a>0$, $x,y\in[0,a]$, $q\geq 0$
\begin{equation}
\int_0^\infty \mathbb{P}_x(X_t\in \D y, \, t< \tau^+_a\wedge \tau^-_0)\D t =\left\{ \frac{W^{(q)} (x) W^{(q)}(a-y)}{W^{(q)}(a)} - W^{(q)}(x-y)\right\}\D y.
\label{SNtwosidedresolve}
\end{equation}
\item[(iii)] Let $a>0$, $x\in[0,a]$, $q\geq 0$ and $f,g$ be positive, bounded measurable functions. Further suppose that $X$ is of bounded variation or $f(0)g(0)=0$.
Then
\begin{multline}
\mathbb{E}_x(e^{-q\tau^-_0} f(X_{\tau^-_0}) g(X_{\tau^-_0 -})\mathbf{1}_{\{\tau^-_0<\tau^+_a\}} ) \\
=\int_0^a\int_{(y,\infty)}f(y-\theta)g(y)
\left\{
\frac{W^{(q)} (x) W^{(q)}(a-y)}{W^{(q)}(a)} - W^{(q)}(x-y)\right\}\Pi(\D\theta)\D y.
\label{two-sided-resolvent}
\end{multline}
\end{itemize}
\end{theorem}


\end{document}